\documentclass{amsart}

\usepackage{amssymb,amscd,amsthm, diagrams}
\usepackage{amscd}
\usepackage{diagrams}

\def\co{\colon\thinspace}

\renewcommand{\phi}{\varphi}
\newcommand{\C}{{\mathbb{C}}}

\newcommand{\Z}{{\mathbb{Z}}}

\newcommand{\ra}{{\rightarrow}}
\newcommand{\lra}{{\longrightarrow}}
\newcommand{\M}{{\widetilde{M}}}

\newtheorem{thm}{Theorem}
\newtheorem{lem}[thm]{Lemma}

\newtheorem{prop}[thm]{Proposition}
      
\newtheoremstyle{note}% name
  {3pt}%      Space above
  {3pt}%      Space below
  {}%         Body font
  {}%         Indent amount (empty = no indent, \parindent = para indent)
  {\scshape}% Thm head font
  {:}%        Punctuation after thm head
  {.5em}%     Space after thm head: " " = normal interword space;
  {}%         Thm head spec (can be left empty, meaning `normal')

\theoremstyle{definition}
\newtheorem*{defn}{Definition}

\theoremstyle{remark}
\newtheorem*{rem}{Remark}
\newtheorem*{rems}{Remarks}
\newtheorem*{nota}{Notation}

\begin{document}
\title[Chern classes of symplectic blow-ups]{A formula for the Chern
classes of\\symplectic blow-ups}

\author{Hansj\"org Geiges and Federica Pasquotto}
\address{Mathematisches Institut, Universit\"at zu K\"oln,
Weyertal 86--90, 50931 K\"oln, Germany}
\email{geiges@math.uni-koeln.de}

\address{Afdeling Wiskunde, Vrije Universiteit, De Boelelaan 1081a,
1081 HV Amsterdam, Netherlands}
\email{pasquott@few.vu.nl}

\date{}

\begin{abstract}
It is shown that the formula for the Chern classes (in the Chow ring)
of blow-ups of algebraic varieties,
due to Porteous and Lascu-Scott, also holds (in the cohomology ring)
for blow-ups of symplectic and complex
manifolds. This was used by the second-named author
in her solution of the geography problem for $8$-dimensional symplectic
manifolds. The proof equally applies to real blow-ups of arbitrary manifolds
and yields the corresponding blow-up formula for the Stiefel-Whitney classes.
In the course of the argument the topological analogue of
Grothendieck's {\em formule clef\/} in intersection theory
is proved.
\end{abstract}

\maketitle

%\setcounter{section}{-1}

%\maketitle

\section{Introduction}
Let $(M, \omega)$ be a symplectic manifold, 
%\marginpar{\small Proposed running head: Chern classes of symplectic blow-ups}
$J$ a tame almost complex structure, that is, $\omega (X,JX)>0$
for any nonvanishing tangent vector~$X$.
The Chern classes of $(M, \omega)$ are defined as the Chern classes of
the complex vector bundle $(TM, J)$.
Since the space of tame almost complex structures for a given symplectic
form is non-empty and contractible (thus in particular connected),
cf.~\cite[p.~65]{mcsa98}, this is a reasonable definition.

Given a symplectic submanifold $N$ of $M$, the normal bundle of $N$
in $M$ carries a complex structure,
and one can then define the blow-up $\M$ of $M$ along~$N$
in analogy with the blow-up of complex manifolds along
complex submanifolds.

The manifold $\M$ admits a symplectic form~$\widetilde{\omega}$,
which coincides with the pullback of $\omega$ outside a small neighbourhood
of the exceptional divisor of the blow-up. The
construction of~$\widetilde{\omega}$, outlined in \cite{gr86} and carried
out in \cite{md84}, depends on a number of choices, which may lead to
non-isomorphic structures. However, the underlying
tame almost complex structures are all homotopic. This allows us to speak
unambiguously of the Chern classes of the symplectic blow-up.

In the algebraic setting, the Chern classes of the blown-up variety
are given by a ``blow-up formula'' found by
Porteous~\cite{port60}. An alternative proof is due
to Lascu and Scott~\cite{lasc73}, cf.\ also \cite{lasc78} and~\cite{lamusc75}.
Here the Chern classes are understood as elements in the Chow ring
(or intersection ring)
of an algebraic variety, see~\cite{fult84} for a brief introduction.
One naturally expects that the blow-up formula
should carry over to the smooth topological
setting, since many formulae in the Chow ring have
analogues in the singular cohomology of manifolds.
All the published proofs, however, depend to some degree
on methods from algebraic geometry that lack an obvious
topological correlate.

In the present paper we provide the necessary translation to
the cohomology of smooth manifolds and use it to show that the blow-up
formula (see Theorem~\ref{thm:blow-up}) equally applies to the blow-up
of symplectic and complex manifolds.
Our proof of the blow-up formula is closest in spirit
to the one in~\cite{lasc78}, but apart from references to the
standard texts \cite{botu82} and \cite{bred93} it is completely
self-contained. We also indicate
how the proof carries over to real blow-ups, where
one obtains the corresponding formula for the Stiefel-Whitney classes.

The proof in~\cite{lasc78} relies in an essential way on
Grothendieck's {\em formule clef\/} in intersection theory
as proved in~\cite{lamusc75}. This is the part where the translation to
the topological setting is least straightforward. Our proof of the
topological analogue of the {\em formule clef\/} uses some ideas
from Quillen's work~\cite{quil71} on complex cobordism theory.
\section{The symplectic blow-up}
\label{section:symp-blow-up}
We briefly recall the definition of the symplectic blow-up;
for details see~\cite{md84}.
Consider a symplectic embedding $i\co (N, \sigma)\ra (M, \omega)$.
We usually identify $N$ with $i(N)\subset M$.
The normal bundle $E$ of $N$ in $M$ 
may be identified with the symplectic orthogonal
bundle of $TN\subset TM$.
Thus $E$ carries a canonical symplectic bundle structure,
given by the restriction of $\omega$ to each fibre, and
hence a homotopically unique tame complex structure as well.
With respect to this structure we
consider the projectivisation $\mathbb{P}(E)$.
Choose a tubular neighbourhood $W$ of $N$ in~$M$. There exists a
closed $2$-form $\rho$ on $E$ which restricts to $\sigma$ along the zero
section and to the canonical symplectic form on each fibre, and with respect
to which $W$ may be symplectically identified with a neighbourhood 
$V$ of the zero section of~$E$. Let $l$ be the tautological line bundle
over $\mathbb{P}(E)$. Denote
by $q$ the bundle projecton $l\ra \mathbb{P}(E)$ and by $\varphi$ the
projection $l\ra E$, so that we have the commutative diagram
\begin{diagram}
l             & \rTo^{\varphi} & E\\
\dTo^q        &                & \dTo_{\pi}\\
\mathbb{P}(E) & \rTo^{p}       & N.  
\end{diagram}
Since $\varphi$ is an isomorphism outside the zero section of $l$,
one can make the following definition, cf.~\cite{md84}.

\begin{defn}
Set $\widetilde{V}:=\varphi^{-1}(V)$; this is a disc sub-bundle
(with fibres real $2$-discs) of the complex line bundle~$l$. The
{\em blow-up} $\widetilde{M}$ of $M$ along $N$ is the manifold
\[ \widetilde{M}:=\overline{M-W}\cup_{\partial \widetilde{V}}
\widetilde{V}.\]
\end{defn}

We shall regard $\mathbb{P}(E)$ as the zero section of
the disc bundle $\widetilde{V}\subset l$, and $N$ as the zero
section of $V\subset E$. The map $\varphi$ gives us
an identification of $\widetilde{V}-\mathbb{P}(E)$
with $V-N$. Thus, we may
alternatively form $\widetilde{M}$ by identifying
$M-N$ and $\widetilde{V}$ along $W-N\cong
V-N\cong\widetilde{V}-\mathbb{P}(E)$.
Either way, we see that
there is a natural inclusion $\widetilde{\imath}\co \mathbb{P}(E)\rightarrow
\M$. This projective space $\mathbb{P}(E)$ is called the
{\em exceptional divisor\/} of the blow-up.

Here is how to construct a symplectic form $\widetilde{\omega}$
on~$\widetilde{M}$. On $M-W$ we set $\widetilde{\omega}=\omega$;
one is then left with defining a symplectic
form on $\widetilde{V}$ which equals
$\varphi^*\rho$ near $\partial\widetilde{V}$.
To do so, one considers a closed
$2$-form $\alpha$ on $\mathbb{P}(E)$ that restricts to the canonical
symplectic form on each fibre of $p$, and that pulls back
under $q^*$ to a form on $l$ that is exact
outside the zero section $\mathbb{P}(E)\subset l$
(such a form may be obtained by the method
of Thurston, cf.\ \cite[Section~6.1]{mcsa98}).
Since $q^*\alpha$ is exact away from the zero section of $l$, one finds a
$1$-form $\beta$ such that $q^*\alpha=d\beta$ on $l-\mathbb{P}(E)$.
There is an $\overline{\varepsilon}>0$, depending on $\rho$ and $\alpha$,
such that for $\varepsilon\in (0,\overline{\varepsilon}]$ 
and with $\lambda$ a radial bump function
on $\widetilde{V}$ which equals $0$ near the boundary, the form
\begin{eqnarray*}
 \tilde{\rho}=
  \left \{ \begin{array}{lll}
  \varphi^*\rho+\varepsilon q^*\alpha & \rm{on} & T\M|_{\mathbb{P}(E)},  \\
  \varphi^*\rho+\varepsilon d(\lambda\beta) & \rm{on} &
     \widetilde{V}-\mathbb{P}(E),
  \end{array} \right.
 \end{eqnarray*}
is nondegenerate on $\widetilde{V}$, and the form
\begin{eqnarray*}
 \widetilde{\omega}=
  \left \{ \begin{array}{lll}
  \omega & \rm{on} & M-W,    \\
  \tilde{\rho} & \rm{on} & \widetilde{V}.
  \end{array} \right.
\end{eqnarray*}
is a symplectic form on $\widetilde{M}$.
\section{The blow-up diagram}
With the symplectic identification of $V\subset E$ with $W\subset M$
understood, we can define a map $f\co\M\rightarrow M$ by
\[
f=
\left \{ \begin{array}{lll}
\rm{id} & \rm{on} & M-W,    \\
\varphi & \rm{on} & \widetilde{V}.
\end{array} \right.
\]
This map is a diffeomorphism outside $\mathbb{P}(E)$, and
$\mathbb{P}(E)=f^{-1}(N)$. In particular, we have
the commutative ``blow-up diagram''
\begin{diagram}
\mathbb{P}(E) & \rTo^{\widetilde{\imath}} & \M\\
\dTo^{p}      &                           & \dTo_{f}\\
N             & \rTo^{i}                  & M.
\end{diagram}

Notice that the normal bundle of $\widetilde{\imath}(\mathbb{P}(E))$ in
$\widetilde{M}$ is isomorphic to $l$. In other words, we have the following
short exact sequence of vector bundles:
\begin{equation}
\label{seq1}
0\lra T\mathbb{P}(E)\lra \widetilde{\imath}^*T\widetilde{M}\lra l\lra 0.
\end{equation}
When there is no ground for confusion, we shall identify
$\mathbb{P}(E)$ with $\widetilde{\imath}(\mathbb{P}(E))\subset\M$.
\section{The Chern classes of symplectic blow-ups}
The construction of a symplectic form on the blow-up of a symplectic
manifold $(M, \omega)$ involves several choices and yields forms which are
not necessarily isomorphic.
Still, we would like to show that the Chern classes of such blown-up
manifolds are well defined.
First we show that we may choose a tame almost complex structure on $M$
that is suitably adapted to the blow-up along~$N$.

\begin{lem}
\label{lem:adaptedJ}
Let $(M, \omega)$ be a symplectic manifold, $N$ a symplectic submanifold
of~$M$. Given any tame almost complex structure $J_0$ on the complement of
$N$ in $M$, one can find a tame almost complex structure $J_M$ on $M$ which
coincides with $J_0$ outside a tubular neighbourhood of $N$ and which is
adapted to $N$, in the sense that $TN$ is $J_M$-invariant and 
the almost complex structure ${J_M}|_{TN}$
is tame with respect to $\omega|_{TN}$.
\end{lem}

\begin{proof}
Let $E$ be the normal bundle of $N$ in $M$ and $W$ a tubular
neighbourhood of $N$, symplectomorphic to a tubular neighbourhood of the
zero section of~$E$. Then $W$ admits an almost complex structure $J_W$
adapted to~$N$.

The space $\mathcal{J}$ of $\omega$-tame almost complex structures on $W-N$
is contractible, i.e.\ the identity map on $\mathcal{J}$ is homotopic to the
constant map sending any almost complex structure $J$ to
$J_0|_{W-N}$. Let
$F\co \mathcal{J}\times I\ra \mathcal{J}$ be the corresponding homotopy.
We may assume that, for some small $\varepsilon >0$, we have
$F(J, t)=J$ for $t\leq\varepsilon$ and $F(J, t)=J_0$
for $t\geq 1-\varepsilon$.

If we let $0\leq t\leq 1$ denote the radial coordinate in $W$, so that
$p\in W$ may be written as $(x, v)$ in some bundle chart, with $x\in N$
and $\| v\| =t$, we
can define an almost complex structure $J_M$ on $M$ as follows:
$$\left \{ \begin{array}{ll}
  J_M(x, v)= F(J_W, t)(x, v) & \quad {\rm for} \quad (x, v)\in W \\ 
  J_M = J_0              & \quad \rm{on} \quad M-W
  \end{array} \right.$$
Then ${J_M}|_{TN}={J_W}|_{TN}$ is again a tame almost complex structure
for $\omega|_{TN}$, and
$J_M=J_0$ outside~$W$.
\end{proof}

By construction, cf.\ \cite[Lemma~3.3]{md84},
the inclusion $\widetilde{\imath}\co\mathbb{P}(E)\ra \M$ is symplectic.
With respect to a tame almost complex structure on $(\M ,\widetilde{\omega})$
adapted to this symplectic submanifold as in Lemma~\ref{lem:adaptedJ},
the sequence~(\ref{seq1}) can be read as an exact sequence of
complex vector bundles. Moreover, the almost complex structures on
$\M$ and $M$ can be chosen in such a way that $f$ is a
pseudoholomorphic map.

We now want to show that such a tame almost complex structure
does not depend, up to homotopy, on the choices in the construction of
a symplectic form $\widetilde{\omega}$ on the blow-up~$\M$.

First of
all, the definition of $\M$ does not depend on the choice of (tame)
complex bundle structure on the normal bundle $E$ of $N$ in~$M$. This
follows from all such choices being homotopic and general bundle
theory, cf.~\cite[Section~4.9]{huse94}.

Now, the construction of $\widetilde{\omega}$ involved the choice
of $2$-forms $\rho$ and $\alpha$, a bump function $\lambda$,
and an $\varepsilon >0$ in an allowable range $(0,\overline{\varepsilon}]$
depending on $\rho$ and~$\alpha$. The conditions on $\rho$ and $\alpha$
are convex. Thus, given two such choices $\rho_i,\alpha_i$, $i=0,1$
(and a corresponding~$\beta_i$), as well as bump functions $\lambda_i$,
one can define $\rho_t,\alpha_t,\beta_t,\lambda_t$, $t\in [0,1]$,
as the respective convex linear combinations $\rho_t=(1-t)\rho_0+t\rho_1$
etc. Since the nondegeneracy condition on $\tilde{\rho}_t$ is
an open condition, and the parameter space $[0,1]$ is compact, we can
find an $\varepsilon >0$ (and smaller than $\overline{\varepsilon}_0,
\overline{\varepsilon}_1$) such that $\tilde{\rho}_t$ is symplectic
for all $t\in [0,1]$. Moreover, varying $\varepsilon$ in the
allowable range $(0,\overline{\varepsilon}]$ gives likewise a family
of symplectic forms. Thus, the corresponding symplectic forms
$\widetilde{\omega}_0$ and $\widetilde{\omega}_1$ on $\M$ are homotopic
through (noncohomologous) symplectic forms, and therefore induce
homotopic tame almost complex structures.

Notice that we could choose an $\varepsilon >0$ only because of our
restriction to a compact family (the convex linear interpolation
between two choices). It is not clear that there is an $\varepsilon >0$
such that one can interpolate between all $\widetilde{\omega}$
through forms corresponding to the same $\varepsilon$.
If this were possible, the interpolation would be through
cohomologous forms, and thus by Moser stability all `$\varepsilon$-blow-ups' 
would be symplectomorphic. This is not known, in general, cf.\
\cite[p.~250]{md84}.
\section{Some cohomological lemmas}
We start by proving some general results, which apply in particular to the
blow-up situation.
Recall that for any map $f\co N\ra M$ of smooth, compact, oriented
manifolds of dimension $n$ and $m$, respectively,
one can define a ``shriek'' or ``transfer'' homomorphism
\[ f^!\co H^{n-p}(N, \partial N) \longrightarrow H^{m-p}(M, \partial M)\]
by $f^!=PD_M\circ f_*\circ PD_N^{-1}$.
Here $PD$
denotes the Poincar\'e duality isomorphism from homology to cohomology.
Likewise, one can define the shriek homomorphism on absolute
cohomology groups, provided $f$ takes the boundary $\partial N$
into~$\partial M$. There is an analogous shriek map $f_!$
on homology, but this will not be used in the present paper.
We shall frequently apply the so-called
{\em projection formula}
\[ f^!(f^*(a)\cup b)=a\cup f^!(b)\]
for $a\in H^*(M)$ or $H^*(M,\partial M)$ and $b\in H^*(N,\partial N)$.
For more details about this and the
following statements see \cite[Sections VI.11, 12 and 14]{bred93}.

\begin{nota}
We use the topologist's convention to label the shriek map
on cohomology with a superscript. More algebraically minded
people sometimes write it with a subscript, emphasising
the covariance of this map. To make matters worse, the corresponding
map on Chow rings is written as $f_*$, whereas the (upper and lower)
shriek maps on Chow rings have a different meaning altogether,
cf.~\cite{fult84}.
\end{nota}

If $W$ is a $k$-disc bundle over a manifold $N$ of dimension $n$,
with projection $\pi\co W\ra N$, and with $i_0\co N\ra W$ the inclusion
of $N$ in $W$ as the zero section, the Thom class of $W$ is defined by
\[ \tau=i_0^!(1)=PD_W\,{i_0}_*[N]\in H^k(W,\partial W). \]
The Thom isomorphism theorem states
that $i_0^!$ is an isomorphism that coincides with the composition
\[ i_0^!\co
H^p(N)\stackrel{\pi^*}{\lra}H^p(W)\stackrel{\cup\tau}{\lra}H^{p+k}(W,
\partial W). \]

If, more generally, $i\co N\ra M$ is a smooth codimension $k$ embedding of
manifolds,
possibly with boundaries (in which case $N$ is required to
meet $\partial M$ transversely in $\partial N$), the Thom class
of the inclusion is defined to be
\[ \tau_N^M=i^!(1)=PD_Mi_*[N]\in H^k(M).\]
Its pull-back under the inclusion $i$ is the Euler class of the normal bundle:
$\chi^M_N=i^*\tau_N^M\in H^k(N)$.
If we denote by $W$ a closed tubular neighbourhood of $N$ in $M$ and identify
it with a $k$-disc sub-bundle of the normal bundle of $N$ in $M$, the Thom
class $\tau_N^M$ is the image of the Thom class $\tau$ of $W$ under
the composition of homomorphisms
\[ H^k(W, \partial W)\stackrel{{\rm exc}^{-1}}{\lra}H^k(M, M-N)\lra H^k(M).
\]
Here exc denotes the excision isomorphism induced by the
inclusion of pairs
\[ (W,\partial W)\ra (M,M-N).\]
(On the level of (co-)homology,
we do not need to distinguish between the pairs $(W,\partial W)$ and
$(W,W-N)$.) The second homomorphism
is induced by the inclusion of $(M,\emptyset )$ in $(M,M-N)$.

From now on, $M$ and $N$ will always be closed manifolds. We regard
$N$ as a submanifold of $M$ and interpret the embedding $i\co N\rightarrow M$
as an inclusion map. The name of the following lemma (and several other
results below) derives from the corresponding statement in intersection
theory.

\begin{lem}[Excision Lemma]
\label{lem:excision}
Let $i\co N\ra M$ be an inclusion of smooth closed manifolds,
$i_c\co M-N\ra M$ the inclusion of the complement of $N$. Suppose
$\lambda\in H^*(M)$ satisfies $i_c^*(\lambda)=0$. Then there exists a class
$\beta\in H^*(N)$ such that $i^!(\beta ) =\lambda$.
\end{lem}

\begin{proof}
As before we write
$i_0\co N\ra W\subset M$ for the inclusion of $N$ in a tubular
neighbourhood~$W$. Write
\[j_M\co (M,\emptyset )\lra (M,M-N) \]
for the inclusion of pairs.
Consider the diagram
\begin{diagram}
H^*(N) & \rTo^{i^!}                            & H^*(M)      & \rTo^{i_c^*} &
                            H^*(M-N) \\
       & \rdTo_{{\rm exc}^{-1} i_0^!} & \uTo>{j_M^*} &            &
                                     \\
       &                                       & H^*(M,M-N). &            &
                                     \\
\end{diagram}

Since the sequence
\[ H^*(M,M-N)\stackrel{j_M^*}{\lra}H^*(M)\stackrel{i_c^*}{\lra}H^*(M-N)\]
is exact, and ${\rm exc}^{-1} i_0^!$ is an isomorphism,
it suffices to show that the diagram above is commutative, i.e.\
$j_M^*{\rm exc}^{-1} i_0^!=i^!$. With the inclusion $i_1\co W\ra M$
we have $i=i_1i_0$, hence $i^!=i_1^!i_0^!$. Since $i_0^!$
is an isomorphism, what we want to show is
\begin{equation}
\label{eqn:excision}
j_M^*{\rm exc}^{-1}=i_1^!,
\end{equation}
in other words, that the following diagram is commutative:
\begin{diagram}
H^*(W, \partial W) & \rTo^{{\rm exc}^{-1}} & H^*(M,M-N) &
       \rTo^{j_M^*}   & H^*(M)  \\
\dTo^{PD_W^{-1}}   &                       &            &
                      & \dTo_{\equiv} \\
H_*(W)             & \rTo^{{i_1}_*}        & H_*(M)     &
       \rTo^{PD_M}    & H^*(M).
\end{diagram}

Since
$H^*(W,\partial W)\cong H^*(W,W-N)$,
a given cohomology class $w\in H^*(W,\partial W)$ can be represented
by a cochain on $W$ that vanishes on singular simplices contained in
$W-N$. Hence we can write $w={\rm exc}(\widetilde{w})$ with
$\widetilde{w}\in H^*(M,M-N)$ represented by a cochain that
vanishes on singular simplices contained in $M-N$.
Notice that if we write, by slight abuse of notation, $i_1$ also
for the inclusion of pairs
\[ i_1\co (W,\partial W)\lra (M,M-N),\]
then
${\rm exc}=i_1^*$. In the following calculations we use $i_1$ in
both senses.

By what we just said about $\widetilde{w}$, we have
\[ \widetilde{w}\cap {i_1}_*[W]=j_M^*(\widetilde{w})\cap [M],\]
since we may represent the fundamental classes $[W]\in H_m(W,\partial W)$
and $[M]\in H_m(M)$ by singular chains that differ only by
singular simplices contained in $M-W$. Hence
\begin{eqnarray*}
PD_M{i_1}_*PD_W^{-1}(w) & = & PD_M{i_1}_*(w\cap [W])\\
        & = & PD_M{i_1}_*(i_1^*(\widetilde{w})\cap [W])\\
        & = & PD_M(\widetilde{w}\cap {i_1}_*[W])\\
        & = & PD_M(j_M^*(\widetilde{w})\cap [M])\\
        & = & j_M^*(\widetilde{w})\\
        & = & j_M^*{\rm exc}^{-1}(w).
\end{eqnarray*}
This concludes the proof.
\end{proof}

\begin{rem}
Equation (\ref{eqn:excision}) explains the statement about the
relation between $\tau_N^M$ and $\tau$ made before the excision lemma:
\[ j_M^*{\rm exc}^{-1}(\tau ) = i_1^!i_0^!(1)=i^!(1)=\tau_N^M.\]
\end{rem}

Up to this point, $N$ was an arbitrary submanifold of $M$. From now on, we
only consider the special set-up described in
Section~\ref{section:symp-blow-up}. We write $2r$ for the rank of
the normal bundle $E$ of $N$ in $M$, and
$c_r(E)$ for the top Chern class (or Euler class)
of this bundle. Then
\[ c_r(E)=\chi^M_N=i^*\tau^M_N=i^*i^!(1).\]
The following lemma generalises this formula.

\begin{lem}[Self-intersection formula]
\label{lem:self}
For any $y\in H^*(N)$ we have
\[ y\cup c_r(E)=i^*i^!(y).\]
\end{lem}

\begin{proof}
Our notation is as in the proof of Lemma~\ref{lem:excision}.
Let $j_W$ be the inclusion of pairs $(W,\emptyset )\ra (W,W-N)$.
Then we have the commutative diagram
\begin{diagram}
H^*(M,M-N)             & \rTo^{j_M^*}  & H^*(M)\\
\dTo^{{\rm exc}=i_1^*} &               & \dTo_{i_1^*} \\
H^*(W,W-N)             & \rTo^{j_W^*}  & H^*(W). 
\end{diagram}

Denote by
\[ \tau\in H^{2r}(W,\partial W)\equiv H^{2r}(W,W-N)\]
the Thom class
$PD_W{i_0}_*[N]$. Notice that $i_0^*j_W^*(\tau )=c_r(E)\in H^{2r}(N)$,
cf.\ \cite[p.~378]{bred93}.

For $w\in H^*(W,\partial W)$ we have, by equation~(\ref{eqn:excision}),
\[ i_1^*i_1^!(w)=i_1^*j_M^*{\rm exc}^{-1}(w)=j_W^*(w).\]
Given $y\in H^*(N)$, we can apply this identity to $w =\pi^*(y)\cup\tau$.
We obtain
\begin{eqnarray*}
i^*i^!(y)  & = &
     i_0^*i_1^*i_1^!i_0^!(y) = i_0^*i_1^*i_1^!(\pi^*(y)\cup\tau )\\
& = & i_0^*j_W^*(\pi^*(y)\cup\tau ) = i_0^*(\pi^*(y)\cup j_W^*(\tau ))\\
& = & y\cup i_0^*j_W^*(\tau ) = y\cup c_r(E).
\end{eqnarray*}
This is the claimed formula.
\end{proof}

In the following lemma (and throughout the remainder this paper) we write
\[ \xi =-c_1(l)\in H^2(\mathbb{P}(E))\]
for the first Chern class of the dual tautological
line bundle~$l^*$; this is the sign convention of~\cite{lasc78}
and motivated by the fact that this $\xi$ is the positive generator
of $H^2(\mathbb{P}(E))$.

\begin{lem}
\label{lem:y-xi}
Suppose $\widetilde{y}\in H^*(\widetilde{M})$ has the property that
$\widetilde{\imath}^*(\widetilde{y})=-\bar{y}\cup\xi$ for some
$\bar{y}\in H^*(\mathbb{P}(E))$.
Then $\widetilde{y}=\widetilde{\imath}^!(\bar{y})+\lambda$ with
$\widetilde{\imath}^*(\lambda )=0$.
\end{lem}

\begin{proof}
We can write $\widetilde{y}=\widetilde{\imath}^!(\bar{y})+(\widetilde{y}
-\widetilde{\imath}^!(\bar{y}))$. By the preceding lemma, applied
to the normal bundle $l$ of $\mathbb{P}(E)\subset\widetilde{M}$,
with top Chern class~$-\xi$, we have
$\widetilde{\imath}^*\widetilde{\imath}^!(\bar{y})=-\bar{y}\cup\xi$. Hence
\[ \widetilde{\imath}^*(\widetilde{y}
-\widetilde{\imath}^!(\bar{y})) = \widetilde{\imath}^*(\widetilde{y})-
\widetilde{\imath}^*\widetilde{\imath}^!(\bar{y})=-\bar{y}\cup\xi
+\bar{y}\cup\xi=0.\; \qed\]
\renewcommand{\qed}{}
\end{proof}

\section{The {\em formule clef}}
The line bundle $l$ may be regarded as a sub-bundle of the pull-back bundle
$p^*E$ over $\mathbb{P}(E)$. The quotient bundle $Q$ is defined by the short
exact sequence
\begin{equation}
\label{seq:Q}
0\lra l \lra p^*E \lra Q \lra 0.
\end{equation}
Recall, cf.~\cite[eqn.~(20.7)]{botu82}, that the cohomology ring
$H^*(\mathbb{P}(E))$ can be described as
\[ H^*(\mathbb{P}(E))=H^*(N)[\xi ]/(\xi^r+p^*c_1(E)\xi^{r-1}+\cdots
+p^*c_r(E)),\]
where from now on we write the cup product as an ordinary product.
The relation
\begin{equation}
\label{eqn:fundamental}
\xi^r+p^*c_1(E)\xi^{r-1}+\cdots +p^*c_r(E)=0
\end{equation}
in $H^*(\mathbb{P}(E))$ will be called the {\em fundamental relation}.
From the exact sequence~(\ref{seq:Q}) we have,
with $c$ denoting the total Chern
class,
\[ p^*c(E)=c(Q)(1-\xi ).\]
By multiplying this equation by $(1+\xi+\cdots +\xi^{r-1})$,
using the fundamental relation, and collecting terms of degree $2r-2$ we find
\begin{equation}
\label{eqn:Q}
c_{r-1}(Q)= \xi^{r-1} +p^*c_1(E)\xi^{r-2}+\cdots +p^*c_{r-1}(E).
\end{equation}

The following key formula, in the algebraic geometric
setting originally conjectured by Grothendieck
and proved in~\cite{lamusc75}, gives an important
tool for computing in the cohomology rings of the manifolds
appearing in the blow-up diagram.

\begin{prop}[Grothendieck's {\em formule clef\/}]
\label{prop:formule-clef}
For any class $y\in H^*(N)$ we have
\[ f^*i^!(y)=\widetilde{\imath}^!(p^*(y)c_{r-1}(Q)).\] 
\end{prop}

The proof of this proposition will take up the rest of this section.
Consider the following commutative diagram, where we
write $D({\mathcal E})$ for the disc bundle associated with
a vector bundle ${\mathcal E}$.
The disc bundles $D(E), D(l)$ will be identified with tubular
neighbourhoods of $N,{\mathbb P}(E)$ in $M,\M$, respectively.
(In other words, $D(E)=W$ and $D(l)=\widetilde{V}$ in our previous
notation.)
The maps in this diagram that have not yet been defined
will be explained presently.
\begin{diagram}
{\mathbb P}(E)     & \rTo^{\widetilde{\imath}_0} & D(l)           &
          \rTo^{\widetilde{\imath}_1} & \M               \\
\dTo^{p_0}         &                             & \dTo_{f_0}     &
                                      & \dTo_{f}         \\
D(T{\mathbb P}(E)) & \rTo^{\overline{\imath}_0}  & D(T{\mathbb P}(E))\oplus
   D(l\oplus Q) &                                    &   \\
\dTo^{p_1}      &                             & \dTo_{f_1}     &
                                      &                  \\
N               & \rTo^{i_0}                  & D(E)              &
          \rTo^{i_1}                  & M                \\ 
\end{diagram}

Here $p_0$ is the inclusion of ${\mathbb P}(E)$ as the zero section in
its tangent disc bundle; $p_1$ is the natural projection of
$D(T{\mathbb P}(E))$
onto ${\mathbb P}(E)$, followed by $p\co {\mathbb P}(E)\ra N$. This means
that $p$ factors as $p=p_1\circ p_0$.

With $D(T{\mathbb P}(E))\oplus D(l\oplus Q)$ we denote the bundle over
${\mathbb P}(E)$ whose fibre over a point $x\in {\mathbb P}(E)$ is the
product of the unit disc in $T_x{\mathbb P}(E)$ with that in
$(l\oplus Q)_x=(p^*E)_x=E_{p(x)}$.
The maps $f_0$ and $f_1$ are the obvious ones.
The restriction of $f$ to $D(l)$ factorises
as $f_1\circ f_0$, so that the square on the
right is indeed commutative.

The square on the top left is commutative since all the maps in that
square are inclusion maps. To see the commutativity of the square on the
bottom left, recall that
\[ p^*E=\{ (x,v)\in {\mathbb P}(E)\times E\co p(x)=\pi (v)\} ,\]
where $\pi\co E\ra N$ denotes the bundle
projection as before. Then for 
$x\in {\mathbb P}(E)$ and $t\in D(T_x{\mathbb P}(E))$ we have
\[ i_0p_1(t)=i_0p(x)=0\in E_{p(x)},\]
and likewise
\[ f_1\overline{\imath}_0(t)=f_1(t,x,0\in E_{p(x)})=0\in E_{p(x)}.\]

In the following lemma and its proof we write $f$ not only for the
map $\M\ra M$, but also for its restriction to subspaces or pairs of
subspaces.

\begin{lem}
\label{lem:fc1}
For $w\in H^*(D(E),\partial D(E))$ we have
\[ f^*i_1^!(w)=\widetilde{\imath}_1^!f^*(w).\]
\end{lem}

\begin{proof}
We have $i_1^!=j_M^*{\rm exc}^{-1}$, see
equation~(\ref{eqn:excision}) in the proof of the excision lemma.
Likewise, with $\widetilde{\jmath}_{\M}$ denoting the
inclusion of pairs
\[ \widetilde{\jmath}_{\M}\co (\M ,\emptyset )\lra (\M ,\M -{\mathbb P}(E)),\]
we have $\widetilde{\imath}_1^!=\widetilde{\jmath}_{\M}^*{\rm exc}^{-1}$,
where ${\rm exc}$ now stands for the excision isomorphism
\[ H^*(\M ,\M -{\mathbb P}(E))\lra H^*(D(l),\partial D(l)).\]
These excision isomorphisms, being induced by inclusions, commute with
$f^*$, and so do $j_M^*$ and $\widetilde{\jmath}_{\M}^*$. This proves
the lemma.
\end{proof}

Next we deal with the square on the bottom left. This is in fact a cartesian
square, i.e.\ $D(T{\mathbb P}(E))$ may be regarded as the fibre product
(or pull-back)
\begin{eqnarray*}
\lefteqn{N\times_{D(E)}(D(T{\mathbb P}(E))\oplus D(p^*E)):=}\\
& & \{ (n;t,x,v)\in N\times D(T{\mathbb P}(E))\times {\mathbb P}(E)\times D(E)
\co\\
& &  t\in T_x{\mathbb P}(E),\, p(x)=\pi (v),\, i_0(n)=f_1(t,x,v)\} .
\end{eqnarray*}
Indeed, we have $f_1(t,x,v)=v$, so the defining equation for the
fibre product becomes $i_0(n)=v$, which implies $v=0\in E_n$. From
$p(x)=\pi (v)$ we then get $n=p(x)$. So the isomorphism of the fibre
product $N\times_{D(E)}(D(T{\mathbb P}(E))\oplus D(p^*E))$
with the disc bundle $D(T{\mathbb P}(E))$ is given by
\[ (n=p(x);t\in T_x{\mathbb P}(E),x,0\in E_n)\longmapsto t\in
T_x{\mathbb P}(E),\]
which has an obvious inverse.

The crucial point for us, however, is the transversality of the maps
$i_0$ and $f_1$. Recall that the Thom class of a disc bundle ${\mathcal D}$
over a closed manifold
is the class in $H^k({\mathcal D},\partial{\mathcal D})$,
with $k$ denoting the fibre dimension,
characterised by the fact that it restricts on each fibre ${\mathcal D}_x$
to the positive generator of $H^k({\mathcal D}_x,\partial {\mathcal D}_x)$.
(All our bundles are complex and thus carry natural
fibre orientations.) 

In our situation we are dealing with a disc bundle
${\mathcal D}:=D(T{\mathbb P}(E))\oplus D(p^*E)$ over the manifold
$D(T{\mathbb P}(E))$, which itself has boundary. So ${\mathcal D}$
is a manifold with corners, but it is still possible to define
a Thom class in this setting:

Write the boundary of ${\mathcal D}$ as $\partial
{\mathcal D}=\partial_B\cup\partial_F$, where the subscripts $B,F$
denote the part of the boundary corresponding to the boundary of
the base and fibre, respectively. The intersection 
$\partial_B\cap\partial_F$ is the codimension~$2$ `corner' of~${\mathcal D}$.
The cap product with the fundamental class of ${\mathcal D}$ gives
a duality isomorphism $PD^{-1}$ from $H^*({\mathcal D},\partial_F)$
to $H_*({\mathcal D},\partial_B)$, see~\cite[p.~358]{bred93}.
So we may define the Thom class of ${\mathcal D}$ as before,
but now this is a cohomology class in $H^k({\mathcal D},\partial_F)$.

The statements about the Thom isomorphism and the characterisation of
the Thom class remain valid with the obvious changes. Thus,
if we write $\overline{\pi}$ for the bundle projection
\[ \overline{\pi}\co {\mathcal D}:=D(T{\mathbb P}(E))\oplus D(p^*E)\lra
D(T{\mathbb P}(E))=:B\]
and $\overline{\tau}\in H^{2r}({\mathcal D},\partial_F)$ for
the Thom class of this disc bundle, the composition
\[ H^p(B)\stackrel{\overline{\pi}^*}{\lra}
H^p({\mathcal D})\stackrel{\cup\overline{\tau}}{\lra}
H^{p+2r}({\mathcal D},\partial_F)\]
is an isomorphism that coincides with
\[  \overline{\imath}_0^!\co H^p(B)
\stackrel{PD^{-1}}{\lra} H_{b-p}(B,\partial B)
\stackrel{i_{0*}}{\lra}H_{b-p}({\mathcal D},\partial_B)
\stackrel{PD}{\lra} H^{p+2r}({\mathcal D},\partial_F),\]
where we wrote $b$ for the dimension of the base.

The characterisation of the Thom class $\overline{\tau}$ as the class that
restricts to the appropriate cohomology generator on each fibre
likewise remains valid --- argue as in the case where the base is a closed
manifold. In our situation this implies $f_1^*(\tau ) =\overline{\tau}$,
with $f_1$ regarded as the map $({\mathcal D},\partial_F)\ra (D(E),\partial
D(E))$.

\begin{lem}[Pull-back]
We have $f_1^*i_0^!=\overline{\imath}_0^!p_1^*$ as
homomorphisms from $H^*(N)$ to $H^*({\mathcal D},\partial_F)$.
\end{lem}

\begin{proof}
From the definitions it is obvious that $\pi\circ f_1=p_1\circ
\overline{\pi}$. Hence, for any class $y\in H^*(N)$ and with
$f_1^*(\tau )=\overline{\tau}$ we get
\[ f_1^*i_0^!(y) = f_1^*(\pi^*(y)\cup\tau )=
\overline{\pi}^*p_1^*(y)\cup\overline{\tau}
=\overline{\imath}_0^!p_1^*(y).\qed\]
\renewcommand{\qed}{}
\end{proof}

Finally, we turn to the square on the top left. For the purposes of our
cohomological computations we may replace the bundle
$D(T{\mathbb P}(E))\oplus D(l\oplus Q)$ by
$D(T{\mathbb P}(E))\oplus D(l)\oplus D(Q)$. Then the maps
$\overline{\imath}_0$ and $f_0$ can be factorised as
\[ \overline{\imath}_0\co D(T{\mathbb P}(E))\stackrel{k_1}{\lra}
D(T{\mathbb P}(E))\oplus D(l)\stackrel{k}{\lra}
D(T{\mathbb P}(E))\oplus D(l)\oplus D(Q)\]
and
\[ f_0\co D(l)\stackrel{k_2}{\lra}
D(T{\mathbb P}(E))\oplus D(l)\stackrel{k}{\lra}
D(T{\mathbb P}(E))\oplus D(l)\oplus D(Q).\]

The commutative diagram
\begin{diagram}
\mathbb{P}(E)      & \rTo^{\widetilde{\imath}_0} & D(l)\\
\dTo^{p_0}         &                             & \dTo_{k_2}\\
D(T{\mathbb P}(E)) & \rTo^{k_1}                  & D(T{\mathbb P}(E))
                                                     \oplus D(l)
\end{diagram}
again constitutes a cartesian square, and the maps $k_1$ and $k_2$
are transverse to each other. By the analogue of the
preceding lemma, we have
\[ k_2^*k_1^!=\widetilde{\imath}_0^!p_0^*.\]

The following lemma and its proof are analogous to an argument
employed by Quillen~\cite{quil71} in his study of the complex cobordism
ring.

\begin{lem}[Clean intersection formula]
For any class $\overline{y}\in H^*(D(T{\mathbb P}(E)))$ we have
\[ f_0^*\overline{\imath}_0^!(\overline{y})=\widetilde{\imath}_0^!(p_0^*
(\overline{y})c_{r-1}(Q)).\]
\end{lem}

\begin{proof}
Write $\widetilde{\pi}$ for the bundle projection $D(l)\ra {\mathbb P}(E)$.
Denote by $\nu_k$ the normal bundle of $D(T{\mathbb P}(E))\oplus D(l)$
in $D(T{\mathbb P}(E))\oplus D(l)\oplus D(Q)$. Then
$k_2^*\nu_k=\widetilde{\pi}^*Q$.

We then compute
\begin{eqnarray*}
f_0^*\overline{\imath}_0^!(\overline{y})
  & = & k_2^*k^*k^!k_1^!(\overline{y})\\
  & = & k_2^*(k_1^!(\overline{y})c_{r-1}(\nu_k))
        \mbox{\hspace{1.98cm}\rm (self-intersection)}\\
  & = & k_2^*k_1^!(\overline{y})k_2^*c_{r-1}(\nu_k)\\
  & = & k_2^*k_1^!(\overline{y})\widetilde{\pi}^*c_{r-1}(Q)\\
  & = & \widetilde{\imath}_0^!p_0^*(\overline{y})\widetilde{\pi}^*c_{r-1}(Q)
       \mbox{\hspace{2.05cm}\rm (pull-back)}\\
  & = & \widetilde{\imath}_0^!(p_0^*(\overline{y})\widetilde{\imath}_0^*
           \widetilde{\pi}^*c_{r-1}(Q))
         \mbox{\hspace{1.5cm}\rm (projection formula)}\\
  & = & \widetilde{\imath}_0^!(p_0^*(\overline{y})c_{r-1}(Q)),
\end{eqnarray*}
which is the desired formula.
\end{proof}

It is now a simple matter to prove the {\em formule clef}.

\begin{proof}[Proof of Proposition~\ref{prop:formule-clef}]
We have 
\begin{eqnarray*}
f^*i^!(y) & = & f^*i_1^!i_0^!(y)\\
          & = & \widetilde{\imath}_1^!f^*i_0^!(y) \mbox{\hspace{3.82cm}
                                   \rm (Lemma~\ref{lem:fc1})}\\
          & = & \widetilde{\imath}_1^!f_0^*f_1^*i_0^!(y)\\
          & = & \widetilde{\imath}_1^!f_0^*\overline{\imath}_0^!p_1^*(y)
                \mbox{\hspace{3.46cm} \rm (pull-back)}\\
          & = & \widetilde{\imath}_1^!\widetilde{\imath}_0^!
                 (p_0^*p_1^*(y)c_{r-1}(Q))
                \mbox{\hspace{2cm} \rm (clean intersection)}\\
          & = & \widetilde{\imath}^! (p^*(y)c_{r-1}(Q)),
\end{eqnarray*}
as was to be shown.
\end{proof}
\section{The blow-up formula}
Write $c({\mathcal E})=1+c_1({\mathcal E})+c_2({\mathcal E})+\ldots$ for
the total Chern class of a complex vector bundle~${\mathcal E}$. If
${\mathcal E}$ is the tangent bundle $TB$ of some manifold $B$, we write
$c(B)$ instead of $c(TB)$.

The expression
\[ \sum_{i=0}^r p^*c_i(E)(1+\xi )^{r-i}(1-\xi )-p^*c(E)
\in H^*({\mathbb P}(E))\]
is obviously of strictly positive degree in~$\xi$, so it makes sense to
divide this term by~$\xi$.

\begin{thm}[Blow-up formula]
\label{thm:blow-up}
Let $N$ be a closed symplectic (resp.\ complex) submanifold of (real)
codimension $2r$ in a symplectic (resp.\ complex) manifold $M$.
Write $E$ for the
normal bundle of $N$ in $M$ with its natural complex bundle
structure, and $\M$ for the symplectic (resp.\ complex)
blow-up of $M$ along $N$.
Let $\xi$ be the first Chern class of the dual tautological line
bundle over the projectivised bundle~${\mathbb P}(E)$. With
$f,\widetilde{i},p$ the maps from the blow-up diagram in
Section~\ref{section:symp-blow-up} we have
\[ c(\M )-f^*c(M)=-\widetilde{\imath}^!
\left[ p^*c(N)\cdot\frac{1}{\xi}\left(
\sum_{i=0}^r p^*c_i(E)(1+\xi )^{r-i}(1-\xi )-p^*c(E)\right)\right] .\] 
\end{thm}

The idea for proving this theorem is as in~\cite{lasc78}
and goes back to Mumford. We first prove
a weak version of the blow-up formula, which differs from the ultimate
version by a term of the form $f^*i^!(\beta )$ with some class
$\beta\in H^*(N)$.
By applying this weak blow-up formula to the blow-up of
$M\times S^2$ along $N\subset M\subset M\times S^2$, one arrives
at Theorem~\ref{thm:blow-up}.

We start with a lemma that allows us to write certain cohomology
classes in $H^*(\M )$ in the form $f^*i^!(\beta )$. This is in fact the only
place where we have to rely on the {\em formule clef\/} (except for an
application in Section~\ref{section:special}).

\begin{lem}
\label{lem:blow-up}
If $\gamma\in H^*({\mathbb P}(E))$ is a cohomology class that
satisfies $\gamma\xi =0$ --- which by the self-intersection formula
is equivalent to saying $\widetilde{\imath}^*
\widetilde{\imath}^!(\gamma )=0$ ---, then there is a class
$\beta\in H^*(N)$ such that $\widetilde{\imath}^!(\gamma )=f^*i^!(\beta )$.
\end{lem}

\begin{proof}
By the structure of the cohomology ring $H^*({\mathbb P}(E))$ described
at the beginning of the preceding section, any class $\gamma\in
H^*({\mathbb P}(E))$ can be described uniquely in the form
\[ \gamma =p^*(\beta _1)\xi^{r-1}+\cdots +p^*(\beta_{r-1})\xi+p^*(\beta_r)\]
with suitable $\beta_i\in H^*(N)$.
Use the fundamental relation~(\ref{eqn:fundamental})
in $H^*({\mathbb P}(E))$ to rewrite the equation $\gamma\xi =0$ as
\[ p^*(\beta_2-\beta_1c_1(E))\xi^{r-1}+\cdots +
p^*(\beta_r-\beta_1c_{r-1}(E))\xi-p^*(\beta_1c_r(E))=0.\]
Since $p^*$ is injective, this implies
\[ \beta_2=\beta_1c_1(E),\, \ldots ,\, \beta_r=\beta_1c_{r-1}(E).\]
With equation~(\ref{eqn:Q}) this yields $p^*(\beta_1 )c_{r-1}(Q)=\gamma$.
By applying the {\em formule clef\/} we obtain
$\widetilde{\imath}^!(\gamma )=f^*i^!(\beta_1)$, so
$\beta:=\beta_1$ is the desired class.
\end{proof}

In order to prove the (weak) blow-up formula, we begin by computing
the result of applying $\widetilde{\imath}^*$ to the left-hand side
of the formula in Theorem~\ref{thm:blow-up}. From the exact sequence
of complex vector bundles
\[ 0\lra TN\lra i^*TM\lra E\lra 0\]
we get
\[ \widetilde{\imath}^*f^*c(M)=p^*i^*c(M)=p^*c(N)p^*c(E).\]
Likewise, from the exact sequence~(\ref{seq1}) we have
\[ \widetilde{\imath}^*c(\M )=c({\mathbb P}(E))c(l)=c({\mathbb P}(E))\cdot
(1-\xi ).\]
Write $V:=\ker (Tp)\subset T{\mathbb P}(E)$ for the bundle of tangent vectors
of ${\mathbb P}(E)$ tangent to the fibres of $p\co {\mathbb P}(E)\ra N$,
so that we have an exact sequence
\begin{equation}
\label{seq:V}
0\lra V\lra T{\mathbb P}(E)\lra p^*TN\lra 0.
\end{equation}
This gives $c({\mathbb P}(E))=c(V)p^*c(N)$.

Moreover, $V$ is isomorphic to the tensor product $Q\otimes l^*$,
cf.~\cite[p.~281]{botu82}. Thus, tensoring the sequence~(\ref{seq:Q})
with $l^*$ yields
\[ 0\lra \C\lra p^*E\otimes l^*\lra V\lra 0;\]
here $\C$ denotes a trivial complex line bundle.
With the formula for computing the total Chern class of the tensor product with
a line bundle \cite[(21.10)]{botu82} we find
\begin{equation}
\label{eqn:V}
c(V)=c(p^*E\otimes l^*)=\sum_{i=0}^rp^*c_i(E)(1+\xi )^{r-i}.
\end{equation}

Putting all this together, we have
\[ \widetilde{\imath}^*(c(\M )-f^*c(M))=
p^*c(N)\left[\sum_{i=0}^r p^*c_i(E)(1+\xi )^{r-i}(1-\xi )-p^*c(E)\right] .\]
Set
\[ \overline{y}= -p^*c(N)\cdot \frac{1}{\xi}
\left[\sum_{i=0}^r p^*c_i(E)(1+\xi )^{r-i}(1-\xi )-p^*c(E)\right]
\in H^*({\mathbb P}(E)) ,\]
so that
\[ \widetilde{\imath}^*(c(\M )-f^*c(M))= -\overline{y}\xi.\]
Lemma~\ref{lem:y-xi} then implies that
\begin{equation}
\label{eqn:blow-up}
c(\M )-f^*c(M)=\widetilde{\imath}^!(\overline{y})+\lambda
\end{equation}
with some class $\lambda\in H^*(\M )$ satisfying $\widetilde{\imath}^*
(\lambda )=0$.

Write $\widetilde{\imath}_c$ for the inclusion $\M -{\mathbb P}(E)
\ra\M$. By the proof of the excision lemma (applied to the
inclusions $\widetilde{\imath}$ and $\widetilde{\imath}_c$) we have
$\widetilde{\imath}_c^*\widetilde{\imath}^!=0$. Hence, by applying
$\widetilde{\imath}_c^*$ to equation~(\ref{eqn:blow-up}) we obtain
\begin{eqnarray*}
\widetilde{\imath}_c^*(\lambda ) & = & \widetilde{\imath}_c^*
                                       (c(\M )-f^*c(M))\\
     & = & c(\M -{\mathbb P}(E))-f^*i_c^*c(M)\\
     & = & c(\M -{\mathbb P}(E))-f^*c(M-N) =0,
\end{eqnarray*}
since $f$ sends $\M -{\mathbb P}(E)$ diffeomorphically
(and pseudoholomorphically) onto $M-N$.
Again by the excision lemma, we know that there is a class
$\gamma\in H^*({\mathbb P}(E))$ with $\widetilde{\imath}^!(\gamma )=\lambda$.
Then $\widetilde{\imath}^*\widetilde{\imath}^!(\gamma)=
\widetilde{\imath}^*(\lambda )=0$, so Lemma~\ref{lem:blow-up}
provides us with a class $\beta\in H^*(N)$ such that
\[ \lambda =\widetilde{\imath}^!(\gamma )=f^*i^!(\beta).\]
Together with equation~(\ref{eqn:blow-up}) this means that we have
proved the formula in Theorem~\ref{thm:blow-up} up to an extra term
$f^*i^!(\beta )$ on the right-hand side. We call this the
{\em weak blow-up formula}.

\vspace{2mm}

We now regard $N$ as a submanifold in $M_S:=M\times S^2$ with its natural
product symplectic structure. The normal bundle of
$N$ in $M_S$ is $E\oplus \C$, with $\C$ denoting a trivial
complex line bundle. Write $\M_S$ for the blow-up of
$M_S$ along $N$, so that we have the following blow-up diagram:
\begin{diagram}
{\mathbb P}(E\oplus\C ) & \rTo^{\widetilde{\imath}_S} & \M_S\\
\dTo^{p_S}              &                             & \dTo_{f_S}\\
N                       & \rTo^{i_S}                  & M_S.
\end{diagram}

Let $l_S$ be the canonical line bundle over ${\mathbb P}(E\oplus\C )$
and set $\xi_S=-c_1(l_S)$. We have
$c_i(E\oplus \C )=c_i(E)$, in particular $c_{r+1}(E\oplus\C )=0$.
So the weak blow-up formula for this set-up reads
\begin{eqnarray}
\label{eqn:weakS}
\lefteqn{c(\M_S )-f_S^*c(M_S)=}\\
& = & -\widetilde{\imath}_S^!
\left[ p_S^*c(N)\cdot\frac{1}{\xi_S}\left(
\sum_{i=0}^r p_S^*c_i(E)(1+\xi_S)^{r+1-i}(1-\xi_S)-p_S^*c(E)\right)\right]
       \nonumber \\
 & & \mbox{}+f_S^*i_S^!(\beta ) \nonumber
\end{eqnarray}
with some class $\beta\in H^*(N)$.

Consider the commutative diagram
\begin{diagram}
\M       & \rTo^{\widetilde{s}} & \M_S\\
\dTo^{f} &                      & \dTo_{f_S}\\
M        & \rTo^s               & M_S,
\end{diagram}
where $s$ and $\widetilde{s}$ are the natural inclusion maps.
We now apply $\widetilde{s}^*$ to the individual summands in
equation~(\ref{eqn:weakS}). In the following computations we write
$\nu_Y^X$ for the normal bundle of a submanifold $Y\subset X$.

We begin with the term $f_S^*i_S^!(\beta )$. Notice that
$i_S=s\circ i$, hence $i_S^!=s^!i^!$. Moreover, by the self-intersection
formula, the composition $s^*s^!$ equals taking the cup product with
$c_1(\nu_{M}^{M_S})=0$. Hence
\begin{equation}
\label{eqn:S1}
\widetilde{s}^*f_S^*i_S^!(\beta )=f^*s^*i_S^!(\beta )=
f^*s^*s^!i^!(\beta )=0.
\end{equation}

Next we deal with the two terms on the left-hand side of~(\ref{eqn:weakS}).
Here we need a lemma.

\begin{lem}
The first Chern class of the normal bundle $\nu_{\M}^{\M_S}$ of
$\M$ in $\M_S$ equals $-\widetilde{\imath}^!(1)$, i.e.\
minus the Thom class of ${\mathbb P}(E)$ in~$\M$.
\end{lem}

\begin{proof}
Let $M'$ be a parallel copy of $M$ in $M_S$, so that $f^{-1}(M')$
is a diffeomorphic copy of $M'$ in $\M_S$. From the explicit
construction of the blow-up one sees that there is a singular chain of
(real) codimension~$1$ in $\M_S$ whose boundary consists of
$\M$ and ${\mathbb P}(E\oplus \C )$ with their natural orientations
(given by the complex structure) and $f^{-1}(M')$ with the reversed
orientation.
This chain can be taken as a smooth
manifold with corner along the transverse intersection
\[ \M\cap {\mathbb P}(E\oplus \C )={\mathbb P}(E).\]

(For the construction of this codimension~$1$ chain,
it is enough to replace $M$ by $E\equiv E\oplus \{ 0\}\subset E\oplus\C$
and $M'$ by a parallel copy $E'\equiv E\oplus\{\varepsilon\}
\subset E\oplus\C$. Then consider the blow-up of
$E\oplus\C$ along the zero section~$N$.
It suffices to deal with
the case where $N$ is a point, where this chain can be
seen quite explicitly. It is best to visualise the blow-up
by cutting out a ball $B^{2r+2}$ centred at zero and of radius smaller
than~$\varepsilon$, and then collapsing its
boundary $S^{2r+1}$ under the Hopf map. A strip $E\times [0,\varepsilon ]$
with boundary $E-E'$
will intersect that $B^{2r+2}$ in half a ball of dimension
$2r+1$. The intersection of $S^{2r+1}$ with $E$ is a $(2r-1)$-dimensional
sphere~$\Sigma$. Collapsing $\Sigma$ gives the blow-up $\M$ of~$M$.
The intersection of $S^{2r+1}$ with $E\times [0,\varepsilon ]$ will be
a $2r$-disc with boundary~$\Sigma$. The interior of
that disc is met by each Hopf fibre of $S^{2r+1}$ exactly once.
That disc will collapse, therefore, to the exceptional divisor
${\mathbb P}(E\oplus\C )$.)

It follows that
\[ [\M ]+[{\mathbb P}(E\oplus \C )]-[f^{-1}(M')]=0\]
in $H_m(\M_S)$, where $m=\dim M$.
Apply $PD_{\M_S}$ to this equation, which gives
\[ \tau_{\M}+\tau_{{\mathbb P}(E\oplus \C )}-\tau_{f^{-1}(M')}=0\]
in $H^2(\M_S)$,
where it is understood that these are the Thom classes of the
respective inclusions into~$\M_S$.

If two submanifolds $A$ and $B$ of a manifold $X$ intersect transversely,
then the pull-back of the Thom class $\tau^X_B$ to $A$ is the
Thom class $\tau^A_{A\cap B}$, see~\cite[pp.~371/2]{bred93}, but
beware the misprint in formula (3) on page 372. Since $\M$ and
${\mathbb P}(E\oplus \C )$ intersect transversely in ${\mathbb P}(E)$,
and $\M$ does not intersect $f^{-1}(M')$, the lemma follows by applying
$\widetilde{s}^*$ to the preceding equation and observing that
$\widetilde{s}^*(\tau_{\M})=\widetilde{s}^*\widetilde{s}^!(1)=
c_1(\nu_{\M}^{\M_S})$ by the self-intersection formula.
\end{proof}

Hence we get
\begin{equation}
\label{eqn:S2}
\widetilde{s}^*c(\M_S)=c(\M)c(\nu_{\M}^{\M_S})=c(\M)\cdot (1-
\widetilde{\imath}^!(1))
\end{equation}
and
\begin{equation}
\label{eqn:S3}
\widetilde{s}^*f_S^*c(M_S)=f^*s^*c(M_S)=f^*(c(M)c(\nu_M^{M_S}))
=f^*c(M).
\end{equation}

Finally, we come to the first summand on the right-hand side
of~(\ref{eqn:weakS}). Consider the following
commutative diagram:
\begin{diagram}
N               & \lTo^{p}    &
{\mathbb P}(E)          & \pile{\rTo^{\widetilde{\imath}_0} \\ \\
                          \lTo_{\widetilde{\pi}}  }      & D(l)             &
      \rTo^{\widetilde{\imath}_1}    & \M \\
\dTo^{\equiv}   &             &
\dTo^{i_{\mathbb{P}}}   &                                & \dTo_{s|_{D(l)}} &
                                     & \dTo_s\\
N               & \lTo^{p_S}  &
{\mathbb P}(E\oplus\C ) & \pile{\rTo^{\widetilde{\imath}_{S0}} \\ \\
                          \lTo_{\widetilde{\pi}_S}  }    & D(l_S)           &
      \rTo^{\widetilde{\imath}_{S1}} & \M_S.
\end{diagram}
Here $i_{\mathbb{P}}$ denotes the natural inclusion of
${\mathbb P}(E)$ in ${\mathbb P}(E\oplus\C )$, and, as before, $l_S$ the
canonical line bundle over ${\mathbb P}(E\oplus\C )$.

We claim that $s^*\widetilde{\imath}_S^!=\widetilde{\imath}^!i_{\mathbb{P}}^*$.
This follows by considering the two squares on the right separately.
Indeed, the equality $s^*\widetilde{\imath}_{S1}^!=
\widetilde{\imath}_1^!(s|_{D(l)})^*$ is proved exactly like
Lemma~\ref{lem:fc1}. The equality $(s|_{D(l)})^*
\widetilde{\imath}_{S0}^!=\widetilde{\imath}_0^!i_{\mathbb{P}}^*$
follows from the observation that the Thom class
$\widetilde{\imath}_{S0}^!(1)$ pulls back under
$(s|_{D(l)})^*$ to the Thom class $\widetilde{\imath}_0^!(1)$;
this in turn is a consequence of $s|_{D(l)}$ being an isomorphism on fibres
and the characterisation of the Thom class as generator of the fibre
(rel boundary) cohomology. Hence
\[ (s|_{D(l)})^*\widetilde{\imath}_{S0}^!(.)=
(s|_{D(l)})^*(\widetilde{\pi}_S^*(.)\widetilde{\imath}_{S0}^!(1))=
\widetilde{\pi}^*i_{\mathbb{P}}^*(.)\widetilde{\imath}_0^!(1)=
\widetilde{\imath}_0^!i_{\mathbb{P}}^*(.).\]
From the two said equalities, the claim is immediate.

Furthermore, we have $i_{\mathbb{P}}^*(\xi_S)=\xi$.
Thus, using all this information when we apply
$\widetilde{s}^*$ to the first summand on the right-hand side of
equation~(\ref{eqn:weakS}), we obtain
\begin{eqnarray}
\lefteqn{-\widetilde{s}^* \widetilde{\imath}_S^!
\left[ p_S^*c(N)\cdot\frac{1}{\xi_S}\left(
\sum_{i=0}^r p_S^*c_i(E)(1+\xi_S)^{r+1-i}(1-\xi_S)-p_S^*c(E)\right)\right]
    =} \\
 & = & 
-\widetilde{\imath}^!i_{\mathbb{P}}^*
\left[ p_S^*c(N)\cdot\frac{1}{\xi_S}\left(
\sum_{i=0}^r p_S^*c_i(E)(1+\xi_S)^{r+1-i}(1-\xi_S)-p_S^*c(E)\right)\right]
   \nonumber \\
& = & -\widetilde{\imath}^!
\left[ p^*c(N)\cdot\frac{1}{\xi}\left(
\sum_{i=0}^r p^*c_i(E)(1+\xi)^{r+1-i}(1-\xi)-p^*c(E)\right)\right] .
   \nonumber 
\end{eqnarray}

This expression equals the right-hand side of the blow-up formula
we are aiming to prove, plus an extra summand
\begin{eqnarray*}
\lefteqn{-\widetilde{\imath}^!\left[ p^*c(N)
\sum_{i=0}^r p^*c_i(E)(1+\xi )^{r-i}(1-\xi )\right] =}\\
& = & -\widetilde{\imath}^!\left[p^*c(N)c(V)c(l)\right]
       \mbox{\hspace{2cm} \rm (equation~(\ref{eqn:V}))}\\
& = & -\widetilde{\imath}^!\left[ c({\mathbb P}(E))c(l)\right]
       \mbox{\hspace{2.6cm} \rm (sequence~(\ref{seq:V}))}\\
& = & -\widetilde{\imath}^!\widetilde{\imath}^*c(\M )
       \mbox{\hspace{3.54cm} \rm (sequence~(\ref{seq1}))}
\end{eqnarray*}

Thus, from the weak blow-up formula and
equations (\ref{eqn:S1}), (\ref{eqn:S2}) and (\ref{eqn:S3})
the blow-up formula follows if we can show that
\[ c(\M )\widetilde{\imath}^!(1)=\widetilde{\imath}^!\widetilde{\imath}^*
c(\M ) .\]
But this is simply the projection formula.

This concludes the proof of Theorem~\ref{thm:blow-up}.

\begin{rems}
(1) This blow-up formula was used by the second author in~\cite{pasq05}
to solve the geography problem for symplectic $8$-manifolds.

(2) Our proof of the blow-up formula carries over to give
the corresponding formula for the Stiefel-Whitney classes of any
real manifold $M$ blown up along a submanifold of codimension~$r$
({\em sic\/}!). Simply read $c_i$ as the Stiefel-Whitney class
$w_i\in H^i(.;{\mathbb Z}_2)$, perform all computations in cohomology
with coefficients in $\Z_2$, and in the final part of the
proof above replace $S^2$ by~$S^1$.
\end{rems}
\section{Special cases} 
\label{section:special}
We now derive explicit expressions for some Chern classes of blow-ups
in a few special cases. We write the formulae so as to allow easy comparison
with the expressions in~\cite[pp.~608--611]{grha94},
where an {\em ad hoc} method is used to derive the corresponding results
for blow-ups of complex manifolds.

\vspace{1mm}

(1) {\em The first Chern class of arbitrary blow-ups:} Since
${\mathbb P}(E)$ is of codimension~$2$ in~$\M$, the homomorphism
$\widetilde{\imath}^!$ increases the cohomological degree by~$2$.
So for the computation of $c_1$ we need to identify the degree~$0$
term inside the square brackets in the blow-up formula. This implies
that $c_1(\M)-f^*c_1(M)$ equals $-\widetilde{\imath}^!(1)$ times
the coefficient of the linear term in $\xi$ in
\[ (1+\xi )^r(1-\xi )=1+(r-1)\xi+\ldots , \]
that is, $-(r-1)\widetilde{\imath}^!(1)$. Notice that
$\widetilde{\imath}^!(1)$ equals the Poincar\'e dual of the class
of the exceptional divisor ${\mathbb P}(E)$ in $\M$, so we have

\[ c_1(\M ) = f^*c_1(M)-(r-1)PD_{\M}[{\mathbb P}(E)].\]

\vspace{1mm}

(2) {\em Blow-up at a point:} Here $N$ is a point and $E$ is trivial, so
(with $\dim M=2r$)
\begin{eqnarray*}
c(\M )-f^*c(M) & = & -\widetilde{\imath}^!\left[ \frac{1}{\xi}
                     \bigl( (1+\xi )^r(1-\xi )-1\bigr)\right] \\
               & = & -\widetilde{\imath}^!
                     \sum_{\nu =0}^{r-1}\left(
                     \left(\begin{array}{c}r\\ \nu +1\end{array}\right) -
                     \left(\begin{array}{c}r\\ \nu\end{array}\right)
                     \right) \xi^{\nu}.
\end{eqnarray*}
Set $\eta = -\widetilde{\imath}^!(1)=PD_{\M}[{\mathbb P}(E)]$. Then,
using the projection and the self-intersection formula, we find
\[ \eta^2 = -\eta\cdot\widetilde{\imath}^!(1)=-\widetilde{\imath}^!
\widetilde{\imath}^*(\eta )=\widetilde{\imath}^!\widetilde{\imath}^*
\widetilde{\imath}^!(1)=-\widetilde{\imath}^!(\xi ),\]
and inductively
\[ \eta^{\nu +1}=-\widetilde{\imath}^!(\xi^{\nu}).\]
Hence
\[ c(\M )-f^*c(M)=\sum_{\nu =1}^r \left(
                     \left(\begin{array}{c}r\\ \nu \end{array}\right) -
                     \left(\begin{array}{c}r\\ \nu -1\end{array}\right)
                     \right) \eta^{\nu}.\]
This is consistent with the well-known fact, cf.~\cite[p.~235]{mcsa98},
that blowing up a point is the same as taking the connected sum
with a copy of $\overline{\C P}^r$, i.e.\ $\C P^r$ with the opposite
of its natural orientation. For instance, we can verify the
formula for the Euler characteristic,
\[ \chi (M\# \overline{\C P}^r)=\chi (M) +\chi (\overline{\C P}^r)-2
=\chi (M) +(r-1),\]
as follows. Since $\xi$ is the positive generator of $H^2({\mathbb P}(E))$,
we have $\eta^r=-\widetilde{\imath}^!(\xi^{r-1})=-PD_{\M}
\widetilde{\imath}_*(1)$,
which is the negative generator of $H^{2r}(\M )$. So the formula
for the Euler characteristic follows from
\[ c_r(\M )-f^*c_r(M) =(1-r)\eta^r.\]

\vspace{1mm}

(3) {\em The second Chern class of a symplectic $6$-manifold blown
up along a $2$-dimensional symplectic submanifold:}
In this case, the blow-up formula becomes
\begin{eqnarray*}
c(\M )-f^*c(M) & = & 
   -\widetilde{\imath}^!
   \left[ p^*c(N)\cdot\frac{1}{\xi}\left(
   \sum_{i=0}^2 p^*c_i(E)(1+\xi )^{2-i}(1-\xi )-p^*c(E)\right)\right] \\
  & = & -\widetilde{\imath}^!
        \left[ p^*c(N)\cdot \left( 1-\xi -(\xi^2+p^*c_1(E)\xi +p^*c_2(E)
        \right) \right] \\
  & = & -\widetilde{\imath}^!\bigl( p^*c(N)\cdot (1-\xi )\bigr)        
   \end{eqnarray*}
by the fundamental relation~(\ref{eqn:fundamental}). Hence
\begin{eqnarray*}
c_2(\M ) & = & f^*c_2(M)-\widetilde{\imath}^!(p^*c_1(N)-\xi )\\
         & = & f^*c_2(M)-\widetilde{\imath}^!(p^*i^*c_1(M)-p^*c_1(E)-\xi )\\
         & = & f^*c_2(M)-\widetilde{\imath}^!\widetilde{\imath}^*f^*c_1(M)
               +\widetilde{\imath}^!c_1(Q)
               \mbox{\hspace{.3cm} \rm (equation (\ref{eqn:Q}))}\\
         & = & f^*c_2(M)-f^*c_1(M)\widetilde{\imath}^!(1)+
               f^*i^!(1)
               \mbox{\hspace{.15cm} \rm (projection formula,
                       {\em formule clef}\/)}\\
         & = & f^*(c_2(M)+PD_M[N])-f^*c_1(M)\cdot PD_{\M}[{\mathbb P}(E)].
\end{eqnarray*}

\end{document}